\documentclass[a4paper,11pt,twoside,reqno]{amsart}

\usepackage[utf8]{inputenc}
\usepackage[plainpages=false,pdfpagelabels=true]{hyperref}
\usepackage{amssymb,amsthm,wasysym}

\newtheorem{Satz}{Theorem}[section]

\newtheorem{Lem}[Satz]{Lemma}

\newtheorem{Cor}[Satz]{Corollary}
\theoremstyle{definition}
\newtheorem{Dfn}[Satz]{Definition}
\newtheorem{Bem}[Satz]{Remark}
\newcommand{\tr}{\operatorname{Tr}}
\newcommand{\hess}{\operatorname{Hess}}

\allowdisplaybreaks[1]

\renewcommand{\epsilon}{\varepsilon}

\newcommand{\R}{\ensuremath{\mathbb{R}}}
\numberwithin{equation}{section}

\usepackage[margin=1in]{geometry}

\title{Some Remarks on Energy inequalities for harmonic maps with potential}
\author{Volker Branding}
\date{\today}
\address{TU Wien\\
Institut für diskrete Mathematik und Geometrie\\
Wiedner Hauptstraße 8–10, A-1040 Wien}
\email[]{volker@geometrie.tuwien.ac.at}
\begin{document}

\begin{abstract}
In this note we discuss how several results characterizing the qualitative behavior
of solutions to the nonlinear Poisson equation can be generalized to harmonic maps
with potential between complete Riemannian manifolds.
This includes gradient estimates, monotonicity formulas and Liouville theorems 
under curvature and energy assumptions.
\end{abstract} 

\maketitle

\section{Introduction and Results}
One of the most studied partial differential equations for a scalar function \(u\colon\R^n\to\R\) is
the \emph{Poisson equation}, that is
\[
\Delta u=f,
\]
where \(f\colon\mathbb{R}\to\mathbb{R}\) is some given function.
If one allows the function \(f\) to also depend on \(u\), that is 
\[
\Delta u=f(u),
\]
one calls the equation \emph{nonlinear Poisson equation}. Following the terminology from the literature
we call a solution \(u\in C^3{(\R^n})\) of the nonlinear Poisson equation \emph{entire}.
For entire solutions of the nonlinear Poisson equation the following results characterizing
its qualitative behavior have been obtained:
\begin{enumerate}
\item Suppose that \(F\in C^2(\R)\) is a nonnegative function that is a potential for \(f\),
that is \(F'(u)=f(u)\). Let \(u\) be an entire, bounded solution of the nonlinear Poisson equation.
Then the following energy inequality holds \cite{MR803255}
\begin{align}
\label{modica-estimate}
|\nabla u|^2\leq 2 F(u).
\end{align}
Such kind of inequalities became known as \emph{Modica-type estimates}.
\item Making use of the Modica-type estimate \eqref{modica-estimate} the following Liouville theorem was given
\cite[Theorem 1]{MR803255}:
Suppose \(F\in C^2(\R)\) is a nonnegative function that is a potential for \(f\)
and \(u\) an entire, bounded solution of the nonlinear Poisson equation.
If \(F(u(x_0))=0\) for some \(x_0\in\R^n\), then \(u\) must be constant.
\item Again, making use of the Modica-type estimate, the following monotonicity formula has been established in
\cite{MR2729079}. Let \(u\colon\R^n\to\R\) be an entire, bounded solution of the nonlinear Poisson equation.
Then the following monotonicity formula holds
\begin{align}
\label{monotonicity-nonlinear-poisson}
\frac{d}{dr}\frac{1}{r^{n-1}}\int_{B_r(x)}(\frac{1}{2}|\nabla u|^2+F(u))dx\geq 0,
\end{align}
where \(B_r(x)\) denotes the ball around the point \(x\in\R^n\) with radius \(r\).
\item Another kind of Liouville theorem was achieved in \cite{MR0289961}: Suppose that \(u\) is an entire solution of
\[
\Delta u=f(u,Du).
\]
If \(\frac{\partial f}{\partial u}\geq 0\) and both \(u\) and \(Du\) are bounded, then \(u\) must be constant.
\item Recently, a maximum principle has been established for solutions of \(\Delta u=\nabla F(u)\) in the
vector valued-case \cite{MR3494327, MR3361722}, that is \(u\colon A\subset\R^n\to\R^m\) , where \(A\subset\R^n\)
is some domain. Here it is assumed that the potential \(F\) vanishes at the boundary of a closed convex set.
\end{enumerate}

In this note we focus on the study of a geometric generalization of the nonlinear Poisson equation,
which leads to the notion of \emph{harmonic maps with potential}.
To this end let \((M,h)\) and \((N,g)\) be two Riemannian manifolds, where we set \(n=\dim M\).
For a smooth map \(\phi\colon M\to N\) we consider the Dirichlet energy of the map,
that is \(\int_M|d\phi|^2dM\). In addition, let \(V\colon N\to\R\) be a smooth scalar function.
We consider the following energy functional
\begin{equation}
\label{energy-functional}
E(\phi)=\int_M\big(\frac{1}{2}|d\phi|^2-V(\phi)\big)dM.
\end{equation}
The Euler-Lagrange equation of the functional \eqref{energy-functional} is given by
\begin{equation}
\label{harmonic-potential}
\tau(\phi)=-\nabla V(\phi),
\end{equation}
where \(\tau(\phi)\in\Gamma(\phi^\ast TN)\) denotes the tension field of the map \(\phi\).
Note that in contrast to the Laplacian acting on functions the tension field of a map
between Riemannian manifolds is a nonlinear operator.
Solutions of \eqref{harmonic-potential} are called \emph{harmonic maps with potential}.
We want to point out that motivated from the physical literature one defines \eqref{energy-functional}
with a minus sign in front of the potential.

Harmonic maps with potential have been introduced in \cite{MR1433176}. 
It is shown that due to the presence of the potential, harmonic maps
with potential can have a qualitative behavior that differs from the one of harmonic maps.
Existence results for harmonic maps with potential have been obtained by the heat flow method \cite{MR1800592}, \cite{MR1979036}
under the assumption that the target has negative curvature.
In addition, an existence result for harmonic maps with potential from compact Riemannian manifolds 
with boundary was obtained in  \cite{MR1680678}, where it is assumed that the image of the map lies inside a convex ball.

Besides the aforementioned existence results there also exist several Liouville theorems for harmonic maps with potential.
For a compact domain manifold \(M\) these were derived by the maximum principle under curvature assumptions in 
\cite[Proposition 4]{MR1433176}.
A Liouville theorem for harmonic maps with potential from a complete noncompact Riemannian manifold and the assumption that the
image of the map \(\phi\) lies inside a geodesic ball is given in \cite{MR1618210}.
A monotonicity formula for harmonic maps with potential together with several Liouville theorems was derived in \cite{MR2929724}.

For functions on Riemannian manifolds several generalizations of the Modica-type estimate 
\eqref{modica-estimate} have been established, see \cite{MR3188740}, \cite{MR1359724}.
These results hold under the assumption that the manifold has positive Ricci curvature.

However, is was also noted that estimates of the form \eqref{modica-estimate} do not hold
if we consider vector-valued functions \cite{MR2083306}, \cite{MR2381198}.

It is the aim of this article to discuss if the results obtained for the nonlinear Poisson equation 
stated in the introduction still hold when considering harmonic maps with potential.

This article is organized as follows: In Section 2 we discuss in which sense the Modica-type estimate
\eqref{modica-estimate} for solutions of the nonlinear Poisson equation can be generalized to
harmonic maps with potential between complete Riemannian manifolds. In the last section we will give 
a Liouville theorem for harmonic maps with potential under curvature and boundedness assumptions.

\section{Energy inequalities for harmonic maps with potential}
Before we turn to deriving energy inequalities let us make the following observation:
\begin{Bem}
If we want to model the trajectory of a point particle in a curved space,
we can make use of harmonic maps with potential from a one-dimensional domain,
which are just geodesics coupled to a potential.
To this end we fix some interval \(I\) and consider a curve \(\gamma\colon I\to N\)
that is a solution of \eqref{harmonic-potential}, which in this case reads
\[
\nabla_{\gamma'}\gamma'=-\nabla V(\gamma).
\]
Here \('\) represents the derivative with respect to the curve parameter, which we will denote by \(s\).
For a curve \(\gamma\) satisfying this equation the total energy is conserved, that is
\[
\frac{1}{2}|\gamma'|^2+V(\gamma)=const.
\]
This can easily be seen by calculating
\[
\frac{d}{ds}(\frac{1}{2}|\gamma'|^2+V(\gamma))=\langle\nabla_{\gamma'}\gamma',\gamma'\rangle+\langle\nabla V(\gamma),\gamma'\rangle=0,
\]
where we used the equation for harmonic maps with potential in the last step. 

This fact is well-known in classical mechanics, that is the mechanics of point particles governed by Newton's law.
The total energy consists of the sum of the kinetic and the potential energy and it 
is conserved when the equations of motion are satisfied. 

However, if the dimension of the domain \(M\) is greater then one, we cannot expect that a statement about
the conservation of the total energy will hold in full generality.
\end{Bem}

We will make of the following Bochner formula for a map \(\phi\colon M\to N\), that is
\begin{align}
\label{bocher-formula}
\Delta \frac{1}{2}|d\phi|^2=|\nabla d\phi|^2+\langle d\phi(\text{Ric}^M(e_i)),d\phi(e_i)\rangle
-\langle R^N(d\phi(e_i),d\phi(e_j))d\phi(e_i),d\phi(e_j)\rangle+\langle\nabla\tau(\phi),d\phi\rangle.
\end{align}
Here \(e_i,i=1,\ldots,n\) is an orthonormal basis of \(TM\).
Throughout this article we make use of the Einstein summation convention, that is we sum over repeated indices.
In addition, by the chain rule for composite maps we find
\begin{align}
\label{bochner-potential}
\Delta V(\phi)=&dV(\tau(\phi))+\hess V(d\phi,d\phi)=-|\nabla V(\phi)|^2+\hess V(d\phi,d\phi),
\end{align}
where we used that \(\phi\) is a solution of \eqref{harmonic-potential} in the second step. 	

In order to obtain the Modica-type estimate \eqref{modica-estimate} for solutions of the scalar nonlinear 
Poisson equation one makes use of the so-called \emph{P-function technique}, which heavily makes use of the maximum principle.
The generalization of the \(P\)-function to harmonic maps with potential is given by
\begin{align}
\label{P-function}
P:=\frac{1}{2}|d\phi|^2+V(\phi).
\end{align}
Unfortunately, it turns out that the P-function does not satisfy a ``nice'' inequality in the case of harmonic maps with potential.
\begin{Lem}
Let \(\phi\colon M\to N\) be a smooth harmonic map with potential.
Then the \(P\)-function \eqref{P-function} satisfies the following inequality
\begin{align}
\label{bochner-p}
\Delta P\geq &\langle d\phi(\text{Ric}^M(e_i)),d\phi(e_i)\rangle-\langle R^N(d\phi(e_i),d\phi(e_j))d\phi(e_i),d\phi(e_j)\rangle \\
\nonumber&+\frac{|\nabla P|^2}{|d\phi|^2}-2\frac{\langle\nabla P,\nabla(V(\phi))}{|d\phi|}+\frac{|\nabla(V(\phi))|^2}{|d\phi|^2}-|\nabla V|^2.
\end{align}
\end{Lem}
\begin{proof}
Using the Bochner-formulas \eqref{bocher-formula}, \eqref{bochner-potential} a direct calculation yields
\begin{align*}
\Delta P=|\nabla d\phi|^2+\langle d\phi(\text{Ric}^M(e_i)),d\phi(e_i)\rangle-\langle R^N(d\phi(e_i),d\phi(e_j))d\phi(e_i),d\phi(e_j)\rangle-|\nabla V|^2.
\end{align*}
In addition, we apply the Kato-inequality and find
\begin{align*}
|\nabla d\phi|^2\geq\big|\nabla |d\phi|\big|^2=\big|\frac{\nabla P-\nabla(V(\phi))}{|d\phi|}\big|^2=\frac{|\nabla P|^2}{|d\phi|^2}-2\frac{\langle\nabla P,\nabla(V(\phi))}{|d\phi|}
+\frac{|\nabla(V(\phi))|^2}{|d\phi|^2}
\end{align*}
yielding the result.
\end{proof}

Let us make some comments about \eqref{bochner-p}:
\begin{Bem}
\begin{enumerate}
 \item If the target has dimension one, then the last two terms on the right hand side in \eqref{bochner-p} cancel each other.
  In this case one can successfully apply the maximum principle under the assumption that the domain has positive Ricci curvature
  giving rise to the Modica-type estimate \eqref{modica-estimate}.
 \item If \(\dim N\geq 2\), then the last two terms in \eqref{bochner-p} will no longer cancel each other.
  Moreover, it is well known by counterexamples, see \cite[Section 2]{MR3427610} and references therein,
  that one cannot expect to obtain a Modica-type estimate in the case that \(\dim N\geq 2\).
\end{enumerate}
\end{Bem}

Since we cannot derive energy inequalities by making use of the techniques that were developed for 
solutions of the scalar nonlinear Poisson equation, we will apply ideas that were used
to derive gradient estimates and Liouville theorems for harmonic maps between complete Riemannian manifolds \cite{MR647905}.
Here, one assumes that the image of the map \(\phi\) lies inside a geodesic ball in the target.

\subsection{Gradient estimates for harmonic maps with potential}
In the following we will make use of the following 
\begin{Lem}
Let \(\phi\colon M\to N\) be a smooth harmonic map with potential.
Suppose that the Ricci curvature of \(M\) and the Hessian of the potential \(V\) satisfy
\(\operatorname{Ric}^M-\hess V\geq-A_V\) and that the sectional curvature \(K^N\) of \(N\) satisfies \(K^N\leq B\).
Then the following inequality holds
\begin{align}
\label{bochner-estimate-dphi}
\frac{\Delta|d\phi|^2}{|d\phi|^2}\geq\frac{1}{2}\frac{\big|d|d\phi|^2\big|^2}{|d\phi|^4}-2A_V-2B|d\phi|^2.
\end{align}
\end{Lem}
\begin{proof}
This follows from the Bochner formula \eqref{bocher-formula} and the identity \(\big|d|d\phi|^2\big|^2\leq 4|d\phi|^2|\nabla d\phi|^2\).
\end{proof}

Now we fix a point \(x_0\) in \(M\) and by \(r\) we denote the Riemannian distance from the point \(x_0\). 
Let \(\eta\colon N\to\R\) be a positive function.
On the geodesic ball \(B_r(x_0)\) in \(M\) we define the function
\begin{equation}
\label{definition-F}
F:=\frac{(a^2-r^2)^2}{(\eta\circ\phi)^2}|d\phi|^2.
\end{equation}
Clearly, the function \(F\) vanishes on the boundary \(B_a(x_0)\),
hence \(F\) attains its maximum at an interior point \(x_{max}\).
We can assume that the Riemannian distance function \(r\) is smooth near the point \(x_{max}\),
see \cite[Section 2]{MR573431}.

In the following we will apply the Laplacian comparison theorem, see \cite[p. 20]{MR521983}, that is
\[
\Delta r^2\leq C_L(1+r)
\]
with some positive constant \(C_L\). Moreover, we make use of the Gauss Lemma, that is \(|dr|^2=1\).

\begin{Lem}
Let \(\phi\colon M\to N\) be a smooth harmonic map with potential.
Suppose that the Ricci curvature of \(M\) and the Hessian of the potential \(V\) satisfy
\(\operatorname{Ric}^M-\hess V\geq-A_V\) and that the sectional curvature \(K^N\) of \(N\) satisfies \(K^N\leq B\).
Then the following inequality holds
\begin{align}
\label{inequality-a}
0\geq&-2A_V-\frac{2C_L(1+r)}{a^2-r^2}-\frac{16r^2}{(a^2-r^2)^2}-\frac{8r|d(\eta\circ\phi)|}{(a^2-r^2)(\eta\circ\phi)}
-2\frac{\Delta(\eta\circ\phi)}{\eta\circ\phi}-2B|d\phi|^2.
\end{align}
\end{Lem}

\begin{proof}
At the maximum \(x_{max}\) the first derivative of \eqref{definition-F} vanishes, yielding
\begin{align}
\label{inequality-b}
0=\frac{-2dr^2}{a^2-r^2}+\frac{d|d\phi|^2}{|d\phi|^2}-\frac{2d(\eta\circ\phi)}{\eta\circ\phi}.
\end{align}
Applying the Laplacian to \eqref{definition-F} at \(x_{max}\) gives
\begin{align}
\label{inequality-c}
0\geq\frac{-2\Delta r^2}{a^2-r^2}-\frac{2|dr^2|^2}{(a^2-r^2)^2}+\frac{\Delta|d\phi|^2}{|d\phi|^2}-\frac{\big|d|d\phi|^2\big|^2}{|d\phi|^4}
-2\frac{\Delta(\eta\circ\phi)}{\eta\circ\phi}+\frac{2|d(\eta\circ\phi)|^2}{(\eta\circ\phi)^2}.
\end{align}
Squaring \eqref{inequality-b} we find
\begin{align}
\label{inequality-d}
\frac{\big|d|d\phi|^2\big|^2}{|d\phi|^4}\leq\frac{4|dr^2|^2}{(a^2-r^2)^2}+\frac{4|d(\eta\circ\phi)|^2}{(\eta\circ\phi)^2}+\frac{8|dr^2||d(\eta\circ\phi)|}{(a^2-r^2)\eta\circ\phi}.
\end{align}
Inserting \eqref{bochner-estimate-dphi} and \eqref{inequality-d} into \eqref{inequality-c} and 
using the Gauss Lemma we get the claim.
\end{proof}

To obtain a gradient estimate from \eqref{inequality-a} for noncompact manifolds \(M\) and \(N\) we have to specify the function \(\eta\).

First, we choose a function \(\eta\) that is adapted to the geometry of the target manifold motivated
by a similar calculation for harmonic maps between complete manifolds \cite{MR647905}.
Let \(\rho\) be the Riemannian distance function from the point \(y_0\) in the target manifold \(N\). 
We define 
\begin{equation}
\label{eta-distance}
\xi:=\sqrt{d}\cos(\sqrt{d}\rho) 
\end{equation}
with some positive number \(\sqrt{d}\) to be fixed later, where \(B_R(y_0)\) denotes the geodesic ball
of radius \(R\) around the point \(y_0\) in \(N\).
We will assume that \(R<\pi/(2\sqrt{d})\), thus \(0<\xi(R)<\sqrt{d}\) on the ball \(B_R(y_0)\). 
\begin{Lem}
On the geodesic ball \(B_R(y_0)\) we have the following estimate
\begin{equation}
\label{hessian-comparison}
\operatorname{Hess}\xi\leq-d^\frac{3}{2}\cos(\sqrt{d}\rho)g,
\end{equation}
where \(g\) denotes the Riemannian metric on \(N\).
\end{Lem}
\begin{proof}
This follows from  the Hessian Comparison theorem, see \cite[Proposition 2.20]{MR521983}
and \cite[p. 93]{MR647905}.
\end{proof}

We will also make use of the following fact:
If \(c_1x^2-c_2x-c_3\leq 0\) for \(c_i>0,i=1,2,3\), then the following inequality holds
\begin{align}
\label{quadratic-inequality}
x\leq\max\{2c_2/{c_1},2\sqrt{c_3/{c_1}}\}. 
\end{align}

At this point we can give the following two results similar to \cite[Theorem 3.2]{MR1680678}:

\begin{Satz}
\label{theorem-ball}
Let \(\phi\colon M\to N\) be a smooth harmonic map with potential.
Suppose that the Ricci curvature of \(M\) and the Hessian of the potential \(V\) satisfy
\(\operatorname{Ric}^M-\hess V\geq-A_V\) and that the sectional curvature \(K^N\) of \(N\) satisfies \(K^N\leq B\).
Moreover, assume that \(\phi(M)\subset B_R(y_0)\), where \(B_R(y_0)\) is the geodesic ball of radius \(0<R<\frac{\pi}{2\sqrt{d}}\) 
around \(y_0\) in \(N\) with \(B<d\).
Then the following estimate holds
\begin{align}
|d\phi|\leq\max\big(\frac{16r\sqrt{d}}{C_2(a^2-r^2)\cos(\sqrt{d}\rho)},
\frac{2}{\sqrt{C_2}}(2A_V+\frac{2C_L(1+r)}{a^2-r^2}+\frac{16r^2}{(a^2-r^2)^2}
+\frac{2\sqrt{d}|\nabla V|}{\cos(\sqrt{d}\rho)})^\frac{1}{2}\big),
\end{align}
where the positive constant \(C_2\) depends on the geometry of \(N\).
\end{Satz}

\begin{proof}
We choose the function \(\xi\) defined in \eqref{eta-distance} and insert it for \(\eta\) in \eqref{inequality-a}.
By the Hessian comparison theorem \eqref{hessian-comparison} we find
\begin{align*}
-\Delta(\xi\circ\phi)=-d\xi(\tau(\phi))-\hess\xi(d\phi,d\phi)\geq-d|\nabla V|+d^\frac{3}{2}\cos(\sqrt{d}\rho)|d\phi|^2,
\end{align*}
where we also used that \(\phi\) is a harmonic map with potential.
In addition, making use of the assumption on the image of \(\phi(M)\),
there exists a positive constant \(C_2\) such that
\begin{align*}
-2B+\frac{2d^\frac{3}{2}\cos(\sqrt{d}\rho)}{\sqrt{d}\cos(\sqrt{d}\rho)}=-2B+2d>C_2
\end{align*}
holds. Inserting this into \eqref{inequality-a} we find
\begin{align*}
0\geq C_2|d\phi|^2-\frac{8rd}{(a^2-r^2)\sqrt{d}\cos(\sqrt{d}\rho)}|d\phi|
-2A_V-\frac{2C_L(1+r)}{a^2-r^2}-\frac{16r^2}{(a^2-r^2)^2}-\frac{2d|\nabla V|}{\sqrt{d}\cos(\sqrt{d}\rho)}.
\end{align*}
The claim then follows from \eqref{quadratic-inequality}.
\end{proof}

\begin{Cor}
Under the assumptions of Theorem \ref{theorem-ball} we can take the limit \(a\to\infty\) while keeping the point \(x_0\) in \(M\)
fixed and obtain the estimate
\begin{align*}
|d\phi|^2\leq C\big(A_V+\frac{\sqrt{d}|\nabla V|}{\cos(\sqrt{d}\rho)}\big).
\end{align*}
If \(M\) has positive Ricci curvature and if the potential \(V(\phi)\) is concave, then 
the following inequality holds
\begin{align*}
\label{estimate-dphi-a}
|d\phi|^2\leq C\frac{\sqrt{d}|\nabla V|}{\cos(\sqrt{d}\rho)},
\end{align*}
which can be interpreted as a Modica-type estimate for harmonic maps with potential.
\end{Cor}

There is another way how we can obtain a gradient estimate from \eqref{inequality-a},
by assuming that the potential \(V(\phi)\) has a special structure. 
More precisely, we have the following 

\begin{Satz}
\label{theorem-energy2}
Let \(\phi\colon M\to N\) be a smooth harmonic map with potential.
Suppose that the Ricci curvature of \(M\) satisfies
\(\operatorname{Ric}^M\geq-A\) and that the sectional curvature \(K^N\) of \(N\) satisfies \(K^N\leq B\).
Moreover, assume that the potential \(V\) satisfies 
\begin{align*}
V(\phi)>0,\qquad -\hess V> BV(\phi) g.
\end{align*}
Then the following energy estimate holds
\begin{align}
|d\phi|\leq\max\big(\frac{16r}{C_3(a^2-r^2)V(\phi)},
\frac{2}{\sqrt{C_3}}(2A+\frac{2C_L(1+r)}{a^2-r^2}+\frac{16r^2}{(a^2-r^2)^2}
\big).
\end{align}
The constant \(C_3\) depends on the geometry of \(N\).
\end{Satz}
\begin{proof}
We make use of the formula \eqref{inequality-a}, where we now choose \(\eta(\phi)=V(\phi)\). 
Making use  of the assumptions on the potential \(V(\phi)\) we note that
\[
-\Delta V(\phi)=-dV(\tau(\phi))-\hess V(d\phi,d\phi)=|\nabla V|^2-\hess V(d\phi,d\phi)>BV(\phi)|d\phi|^2.
\]
In addition, again by the assumptions on the potential \(V(\phi)\), we get
\begin{align*}
-2B|d\phi|^2-2\frac{\hess V(d\phi,d\phi)}{V(\phi)}>C_3|d\phi|^2
\end{align*}
for some positive constant \(C_3\).
Inserting into \eqref{inequality-a} then yields
\begin{align*}
0\geq&C_3|d\phi|^2-\frac{8r|\nabla V|}{(a^2-r^2)V(\phi)}|d\phi|-2A-\frac{2C_L(1+r)}{a^2-r^2}
-\frac{16r^2}{(a^2-r^2)^2}.
\end{align*}
The statement follows from applying \eqref{quadratic-inequality} again.
\end{proof}

\begin{Cor}
Under the assumptions of Theorem \ref{theorem-energy2} we can take the limit \(a\to\infty\) while keeping the point \(x_0\) in \(M\)
fixed and obtain the estimate
\begin{align*}
|d\phi|^2\leq CA.
\end{align*}
If \(M\) has nonnegative Ricci curvature then \(\phi\) is trivial.
\end{Cor}

\subsection{Generalized Monotonicity formulas}
In the following we will make use of the stress-energy-tensor for harmonic maps with potential, which is locally given by
\begin{align}
\label{stress-energy}
S_{ij}=\frac{1}{2}|d\phi|^2h_{ij}-\langle d\phi(e_i),d\phi(e_j)\rangle-V(\phi)h_{ij}.
\end{align}
The stress-energy-tensor is divergence-free, when \(\phi\) is a smooth harmonic map with potential \cite{MR2929724}, that is
\[
\nabla_iS_{ij}=0.
\]
Let us recall the following facts:
A vector field \(X\) is called \emph{conformal} if
\[
\mathcal{L}_Xh=fh,
\]
where \(\mathcal{L}\) denotes the Lie-derivative of the metric \(h\) with respect to \(X\) and
\(f\colon M\to\mathbb{R}\) is a smooth function.
\begin{Lem}
Let \(T\) be a symmetric 2-tensor. For a conformal vector field \(X\) the following formula holds
\begin{equation}
\label{conformal-vf}
\operatorname{div}(\iota_X T)=\iota_X\operatorname{div} T+\frac{1}{n}\operatorname{div}X\tr T.
\end{equation}
\end{Lem}
By integrating over a compact region \(U\) and making use of Stokes theorem, we obtain:
\begin{Lem}
Let \((M,h)\) be a Riemannian manifold and \(U\subset M\) be a compact region with smooth boundary.
Then, for any symmetric \(2\)-tensor and a conformal vector field \(X\) the following formula holds
\begin{equation*}
\label{gauss-tensor-formula}
\int_{\partial U}T(X,\nu)d\sigma=\int_U\iota_X\operatorname{div} Tdx+\int_U\frac{1}{n}\operatorname{div}X\tr T dx,
\end{equation*}
where \(\nu\) denotes the normal to \(U\).
\end{Lem}

We now derive a type of monotonicity formula for smooth solutions of \eqref{harmonic-potential} for the domain being \(\mathbb{R}^n\).
\begin{Lem}
Let \(\phi\colon\R^n\to N\) be a smooth harmonic map with potential.
Let \(B_r(x)\) be a ball with radius \(r\) in \(\R^n\). Then the following formula holds
\begin{align}
\label{pre-mono-rn}
r\int_{\partial B_r(x)}\big(\frac{1}{2}|d\phi|^2-V(\phi))d\sigma
-r\int_{\partial B_r(x)}\big|\frac{\partial\phi}{\partial r}\big|^2d\sigma
=(n-2)\int_{B_r(x)}\big(\frac{1}{2}|d\phi|^2-V(\phi))dx\\
\nonumber-2\int_{B_r(x)}V(\phi)dx.
\end{align}
\end{Lem}
\begin{proof}
For \(M=\mathbb{R}^n\) we choose the conformal vector field \(X=r\frac{\partial}{\partial r}\) with \(r=|x|\).
Note that \(\operatorname{div}X=n\).
The statement then follows from \eqref{gauss-tensor-formula} applied to \eqref{stress-energy}.
\end{proof}

Making use of the coarea formula we obtain the following
\begin{Satz}
Let \(\phi\colon\R^n\to N\) be a smooth harmonic map with potential.
Let \(B_r(x)\) be a ball with radius \(r\) in \(\R^n\). Then the following formula holds
\begin{align*}
\frac{d}{dr}r^{2-n}\int_{B_r(x)}(\frac{1}{2}|d\phi|^2-V(\phi))dx
=r^{2-n}\int_{\partial B_r(x)}\big|\frac{\partial\phi}{\partial r}\big|^2dx
-2r^{1-n}\int_{B_r(x)}V(\phi)dx.
\end{align*}
\end{Satz}

\begin{Cor}
Let \(\phi\colon\R^n\to N\) be a smooth harmonic map with potential.
Suppose that \(V(\phi)\leq 0\). Then we have the following monotonicity formula
\begin{align*}
\frac{d}{dr}\frac{1}{r^{n-2}}\int_{B_r(x)}(\frac{1}{2}|d\phi|^2-V(\phi))dx\geq 0.
\end{align*}
Note that this monotonicity formula is different from the one for solutions of the nonlinear Poisson equation \eqref{monotonicity-nonlinear-poisson}
since we do not have a Modica-type estimate for harmonic maps with potential.
\end{Cor}

Monotonicity formulas for harmonic maps with potential with the domain being a Riemannian manifold
have been established in \cite{MR2929724}.

The results presented above also hold for harmonic maps with potential that have lower regularity.
To this end we need the notion of stationary harmonic maps with potential.
\begin{Dfn}
A weak harmonic map with potential is called \emph{stationary harmonic map with potential}
if it is also a critical point of the energy functional with respect to variations of the metric on 
the domain \(M\), that is
\begin{equation}
\label{stationary-condition}
0=\int_Mk^{ij}(\frac{1}{2}|d\phi|^2h_{ij}-\langle d\phi(e_i),d\phi(e_j)\rangle-V(\phi)h_{ij})dM.
\end{equation}
Here \(k^{ij}\) is a smooth symmetric 2-tensor.
\end{Dfn}

Every smooth harmonic map with potential is stationary, 
which is due to the fact that the associated stress-energy-tensor is conserved.
However, a stationary harmonic map with potential can have lower regularity.

For stationary harmonic maps with potential we have the following result generalizing \cite[Theorem 3.1]{MR2729079}:
\begin{Satz}
Let \(\phi\in W^{1,2}_{loc}(M,N)\cap L^\infty_{loc}(M,N)\) be a harmonic map with potential.
Suppose that \(M=\mathbb{R}^n,\mathbb{H}^n\) with \(\dim M\geq 3\) and  
\begin{align*}
\int_M(|d\phi|^2+|V(\phi)|)dM<\infty,
\end{align*}
then the following inequality holds
\begin{align*}
\int_M|d\phi|^2dM\leq\frac{n}{n-2}\int_MV(\phi)dM.
\end{align*}
In particular, this implies that \(\phi\) is constant when \(V(\phi)\leq 0\).
\end{Satz}
\begin{proof}
We will prove the result for the case that \(M=\R^n\).
Let \(\eta\in C_0^\infty(\mathbb{R})\) be a smooth cut-off function satisfying \(\eta=1\) for \(r\leq R\),
\(\eta=0\) for \(r\geq 2R\) and \(|\eta'(r)|\leq\frac{C}{R}\). In addition, 
we choose \(Y(x):=x\eta(r)\in C_0^\infty(\mathbb{R}^n,\mathbb{R}^n)\) with \(r=|x|\).
Hence, we find
\[
k_{ij}=\frac{\partial Y_i}{\partial x^j}=\delta_{ij}\eta(r)+\frac{x_i x_j}{r}\eta'(r).
\]
Inserting this choice into \eqref{stationary-condition} we obtain
\begin{align*}
\int_{\mathbb{R}^n}\eta(r)(((2-n))|d\phi|^2+2nV(\phi))dM=
\int_{\mathbb{R}^n} r\eta'(r)(|d\phi|^2-2\big|\frac{\partial\phi}{\partial r}\big|^2-2V(\phi))dM.
\end{align*}
We can bound the right-hand side as follows
\[
\int_{\mathbb{R}^n} r\eta'(r)(|d\phi|^2-2\big|\frac{\partial\phi}{\partial r}\big|^2-2V(\phi))dM
\leq C\int_{B_{2R}\setminus B_R}(|d\phi|^2+|V(\phi)|)dx.
\]
Making use of the properties of the cut-off function \(\eta\) we obtain
\[
\int_{B_r}((2-n)|d\phi|^2+2nV(\phi))dx\leq C\int_{B_{2R}\setminus B_R}(|d\phi|^2+|V(\phi)|)dx.
\]
Taking the limit \(R\to\infty\) and making use of the assumptions we find
\[
(2-n)\int_{\mathbb{R}^n}(|d\phi|^2+nV(\phi))dM\leq 0,
\]
which finishes the proof for the case that \(M=\R^n\). Making use of the Theorem of Cartan-Hadamard the proof 
carries over to hyperbolic space.
\end{proof}

\begin{Bem}
The last Theorem can be interpreted as an integral version of \eqref{modica-estimate} for bounded harmonic 
maps with potential.
\end{Bem}

Now we derive a generalized monotonicity formula for harmonic maps with potential,
where we take into account the pointwise gradient estimate \eqref{estimate-dphi-a}.

\begin{Satz}
\label{theorem-mono-imrpoved}
Let \(\phi\colon \R^n\to N\) be a smooth harmonic map with potential.
Suppose that the Hessian of the potential \(V\) satisfies
\(-\hess V\geq-A_V\) and that the sectional curvature \(K^N\) of \(N\) satisfies \(K^N\leq B\).
Moreover, assume that \(\phi(M)\subset B_R(y_0)\), where \(B_R(y_0)\) is the geodesic ball of radius \(0<R<\frac{\pi}{2\sqrt{d}}\) 
around \(y_0\) in \(N\) with \(B<d\).

Then the following monotonicity-type formula holds
\begin{align}
\label{mono-imrpoved}
\frac{d}{dr}r^{-n}\bigg(\int_{B_r(x)}(\frac{1}{2}|d\phi|^2-V(\phi))dx\bigg)\geq 
-Cr^{-n-1}\int_{B_r(x)}\big(A_V+\frac{\sqrt{d}|\nabla V|}{\cos(\sqrt{d}\rho)}\big)dx,
\end{align}
where the positive constant \(C\) depends on \(B\).
\end{Satz}

\begin{proof}
Throughout the proof we set
\begin{align*}
e_V(\phi):=\frac{1}{2}|d\phi|^2-V(\phi).
\end{align*}
Making use of the coarea formula and rewriting \eqref{pre-mono-rn} we find
\begin{align*}
r\frac{d}{dr}\int_{B_r(x)}e_V(\phi)dx=r\int_{\partial B_r(x)}\big|\frac{\partial\phi}{\partial r}\big|^2d\sigma
+n\int_{B_r(x)}e_V(\phi)dx-\int_{B_r(x)}|d\phi|^2dx.
\end{align*}
Applying \eqref{estimate-dphi-a} we obtain the following inequality
\begin{align*}
r\frac{d}{dr}\int_{B_r(x)}e_V(\phi)dx-n\int_{B_r(x)}e_V(\phi)dx\geq-C\int_{B_r(x)}\big(A_V+\frac{\sqrt{d}|\nabla V|}{\cos(\sqrt{d}\rho)}\big)dx,
\end{align*}
from which we get the claim.
\end{proof}

Let us make several comments on Theorem \ref{theorem-mono-imrpoved}:
\begin{Bem}
\begin{enumerate}
\item The monotonicity type-formula \eqref{mono-imrpoved} can be interpreted as the generalization of \eqref{monotonicity-nonlinear-poisson}
to harmonic maps with potential.
\item It is straightforward to generalize \eqref{pre-mono-rn} to the case of the domain being a Riemannian manifold.
\end{enumerate}
\end{Bem}

\section{A Liouville theorem}
In this section we derive a Liouville theorem for harmonic maps with potential from complete 
noncompact manifolds with positive Ricci curvature. Our result is motivated from a similar result for harmonic maps, 
see \cite[Theorem 1]{MR0438388}. In addition, this result also generalizes the Liouville theorem 
for solutions of the nonlinear Poisson equation \cite{MR0289961}, which is stated in detail in the introduction.
\begin{Satz}
\label{theorem-liouville}
Let \((M,h)\) be a complete noncompact Riemannian manifold with nonnegative Ricci curvature,
and \((N,g)\) a manifold with nonpositive sectional curvature. Let \(\phi\colon M\to N\) be a smooth harmonic map with potential
with finite energy, that is \(e(\phi):=\frac{1}{2}|d\phi|^2<\infty\). If the potential \(V\) is concave then \(\phi\) is a constant map. 
\end{Satz}
\begin{proof}
We follow the presentation in \cite[pp. 26]{MR1391729}
For a solution \(\phi\) of \eqref{harmonic-potential} and by the standard Bochner formula \eqref{bocher-formula} we find
\begin{align}
\label{bochner-dphi}
\Delta\frac{1}{2}|d\phi|^2=&|\nabla d\phi|^2+\langle d\phi(\text{Ric}^M(e_i)),d\phi(e_i)\rangle
-\langle R^N(d\phi(e_i),d\phi(e_j))d\phi(e_i),d\phi(e_j)\rangle \\
\nonumber&-\hess V(d\phi,d\phi).
\end{align}
Making use of the curvature assumptions and the fact that the potential is a concave function, \eqref{bochner-dphi} yields
\begin{equation}
\label{liouville-complete-a}
\Delta e(\phi)\geq |\nabla d\phi|^2.
\end{equation}
In addition, by the Cauchy-Schwarz inequality we find
\begin{equation}
\label{liouville-complete-b}
|de(\phi)|^2\leq 2e(\phi)|\nabla d\phi|^2.
\end{equation}
We fix a positive number \(\epsilon>0\) and calculate
\begin{align*}
\Delta\sqrt{e(\phi)+\epsilon}=\frac{\Delta e(\phi)}{2\sqrt{e(\phi)+\epsilon}}
-\frac{1}{4}\frac{|de(\phi)|^2}{(e(\phi)+\epsilon)^\frac{3}{2}} 
\geq\frac{|\nabla d\phi|^2}{2\sqrt{e(\phi)+\epsilon}}\big(1-\frac{e(\phi)}{e(\phi)+\epsilon}\big)
\geq 0,
\end{align*}
where we used \eqref{liouville-complete-a} and \eqref{liouville-complete-b}.
Let \(\eta\) be an arbitrary function on \(M\) with compact support. 
We obtain
\begin{align*}
0\leq&\int_M\eta^2\sqrt{e(\phi)+\epsilon}\Delta\sqrt{e(\phi)+\epsilon}dM \\
=&-2\int_M\eta\sqrt{e(\phi)+\epsilon}\langle\nabla\eta,\nabla\sqrt{e(\phi)+\epsilon}\rangle dM
-\int_M\eta^2|\nabla\sqrt{e(\phi)+\epsilon}|^2dM.
\end{align*}
Now let \(x_0\) be a point in \(M\) and let \(B_R,B_{2R}\) be geodesic balls centered at \(x_0\)
with radii \(R\) and \(2R\).
We choose a cutoff function \(\eta\) satisfying
\[
\eta(x)=
\begin{cases}
1,\qquad x\in B_R,\\
0, \qquad x\in M\setminus B_{2R}.
\end{cases}
\]
In addition, we choose \(\eta\) such that
\[
0\leq\eta\leq 1,\qquad |\nabla\eta|\leq\frac{C}{R}
\]
for a positive constant \(C\).
Then, we find
\begin{align*}
0\leq&-2\int_{B_{2R}}\eta\sqrt{e(\phi)+\epsilon}\langle\nabla\eta,\nabla\sqrt{e(\phi)+\epsilon}\rangle dx
-\int_{B_{2R}}\eta^2|\nabla\sqrt{e(\phi)+\epsilon}|^2 dx\\
\leq&2\big(\int_{B_{2R\setminus B_R}}\eta^2|\sqrt{e(\phi)+\epsilon}|^2dx\big)^\frac{1}{2}
\big(\int_{B_{2R\setminus B_R}}|\nabla\eta|^2(e(\phi)+\epsilon)dx\big)^\frac{1}{2}\\
&-\int_{B_{2R}{\setminus B_R}}\eta^2|\nabla\sqrt{e(\phi)+\epsilon}|^2dx
-\int_{B_R}|\nabla\sqrt{e(\phi)+\epsilon}|^2dx.
\end{align*}
We therefore obtain
\begin{align*}
\int_{B_R}|\nabla\sqrt{e(\phi)+\epsilon}|^2dx\leq\int_{B_{2R}{\setminus B_R}}|\nabla\eta|^2(e(\phi)+\epsilon)dx
\leq\frac{C^2}{R^2}\int_{B_{2R}}(e(\phi)+\epsilon)dx.
\end{align*}
We set \(B'_R:=B_R\setminus\{x\in B_R\mid e(\phi)(x)=0\}\) and find
\begin{align*}
\int_{B'_R}\frac{|\nabla(e(\phi)+\epsilon)|^2}{4(e(\phi)+\epsilon)}dx
\leq\frac{C^2}{R^2}\int_{B_{2R}}(e(\phi)+\epsilon)dx.
\end{align*}
Letting \(\epsilon\to 0\) we get
\begin{align*}
\int_{B'_R}\frac{|\nabla e(\phi)|^2}{4e(\phi)}dx
\leq\frac{C^2}{R^2}\int_{B_{2R}}e(\phi)dx.
\end{align*}
Now, letting \(R\to\infty\) and under the assumption that the energy is finite, we have
\[
\int_{M\setminus\{e(\phi)=0\}}\frac{|\nabla e(\phi)|^2}{4e(\phi)}dM\leq 0,
\]
hence the energy \(e(\phi)\) has to be constant.
If \(e(\phi)\neq 0\), then the volume of \(M\) would have to be finite.
However, by \cite[Theorem 7]{MR0417452} the volume of a complete and noncompact Riemannian 
manifold with nonnegative Ricci curvature is infinite. Hence \(e(\phi)=0\),
which yields the result.
\end{proof}
\bibliographystyle{plain}
\bibliography{mybib}

\def\cprime{$'$}
\begin{thebibliography}{10}

\bibitem{MR2729079}
Nicholas~D. Alikakos.
\newblock Some basic facts on the system {$\Delta u-W_u(u)=0$}.
\newblock {\em Proc. Amer. Math. Soc.}, 139(1):153--162, 2011.

\bibitem{MR3361722}
Nicholas~D. Alikakos and Giorgio Fusco.
\newblock A maximum principle for systems with variational structure and an
  application to standing waves.
\newblock {\em J. Eur. Math. Soc. (JEMS)}, 17(7):1547--1567, 2015.

\bibitem{MR3494327}
Panagiotis Antonopoulos and Panayotis Smyrnelis.
\newblock A maximum principle for the system {$\Delta u-\nabla W(u)=0$}.
\newblock {\em C. R. Math. Acad. Sci. Paris}, 354(6):595--600, 2016.

\bibitem{MR1618210}
Qun Chen.
\newblock Liouville theorem for harmonic maps with potential.
\newblock {\em Manuscripta Math.}, 95(4):507--517, 1998.

\bibitem{MR1680678}
Qun Chen.
\newblock Maximum principles, uniqueness and existence for harmonic maps with
  potential and {L}andau-{L}ifshitz equations.
\newblock {\em Calc. Var. Partial Differential Equations}, 8(2):91--107, 1999.

\bibitem{MR1979036}
Qun Chen and Zhen-Rong Zhou.
\newblock Heat flows of harmonic maps with potential into manifolds with
  nonpositive curvature.
\newblock {\em Arch. Math. (Basel)}, 80(2):216--224, 2003.

\bibitem{MR573431}
Shiu~Yuen Cheng.
\newblock Liouville theorem for harmonic maps.
\newblock In {\em Geometry of the {L}aplace operator ({P}roc. {S}ympos. {P}ure
  {M}ath., {U}niv. {H}awaii, {H}onolulu, {H}awaii, 1979)}, Proc. Sympos. Pure
  Math., XXXVI, pages 147--151. Amer. Math. Soc., Providence, R.I., 1980.

\bibitem{MR647905}
Hyeong~In Choi.
\newblock On the {L}iouville theorem for harmonic maps.
\newblock {\em Proc. Amer. Math. Soc.}, 85(1):91--94, 1982.

\bibitem{MR1433176}
Ali Fardoun and Andrea Ratto.
\newblock Harmonic maps with potential.
\newblock {\em Calc. Var. Partial Differential Equations}, 5(2):183--197, 1997.

\bibitem{MR1800592}
Ali Fardoun, Andrea Ratto, and Rachid Regbaoui.
\newblock On the heat flow for harmonic maps with potential.
\newblock {\em Ann. Global Anal. Geom.}, 18(6):555--567, 2000.

\bibitem{MR2083306}
Alberto Farina.
\newblock Two results on entire solutions of {G}inzburg-{L}andau system in
  higher dimensions.
\newblock {\em J. Funct. Anal.}, 214(2):386--395, 2004.

\bibitem{MR521983}
R.~E. Greene and H.~Wu.
\newblock {\em Function theory on manifolds which possess a pole}, volume 699
  of {\em Lecture Notes in Mathematics}.
\newblock Springer, Berlin, 1979.

\bibitem{MR2381198}
Changfeng Gui.
\newblock Hamiltonian identities for elliptic partial differential equations.
\newblock {\em J. Funct. Anal.}, 254(4):904--933, 2008.

\bibitem{MR2929724}
Hezi Lin, Guilin Yang, Yibin Ren, and Tian Chong.
\newblock Monotonicity formulae and {L}iouville theorems of harmonic maps with
  potential.
\newblock {\em J. Geom. Phys.}, 62(9):1939--1948, 2012.

\bibitem{MR3188740}
Li~Ma and Ingo Witt.
\newblock Liouville theorem for the nonlinear {P}oisson equation on manifolds.
\newblock {\em J. Math. Anal. Appl.}, 416(2):800--804, 2014.

\bibitem{MR803255}
Luciano Modica.
\newblock A gradient bound and a {L}iouville theorem for nonlinear {P}oisson
  equations.
\newblock {\em Comm. Pure Appl. Math.}, 38(5):679--684, 1985.

\bibitem{MR1359724}
Andrea Ratto and Marco Rigoli.
\newblock Gradient bounds for {L}iouville's type theorems for the {P}oisson
  equation on complete {R}iemannian manifolds.
\newblock {\em Tohoku Math. J. (2)}, 47(4):509--519, 1995.

\bibitem{MR0438388}
Richard Schoen and Shing~Tung Yau.
\newblock Harmonic maps and the topology of stable hypersurfaces and manifolds
  with non-negative {R}icci curvature.
\newblock {\em Comment. Math. Helv.}, 51(3):333--341, 1976.

\bibitem{MR0289961}
James Serrin.
\newblock Entire solutions of nonlinear {P}oisson equations.
\newblock {\em Proc. London. Math. Soc. (3)}, 24:348--366, 1972.

\bibitem{MR3427610}
Panayotis Smyrnelis.
\newblock Gradient estimates for semilinear elliptic systems and other related
  results.
\newblock {\em Proc. Roy. Soc. Edinburgh Sect. A}, 145(6):1313--1330, 2015.

\bibitem{MR1391729}
Yuanlong Xin.
\newblock {\em Geometry of harmonic maps}.
\newblock Progress in Nonlinear Differential Equations and their Applications,
  23. Birkh\"auser Boston Inc., Boston, MA, 1996.

\bibitem{MR0417452}
Shing~Tung Yau.
\newblock Some function-theoretic properties of complete {R}iemannian manifold
  and their applications to geometry.
\newblock {\em Indiana Univ. Math. J.}, 25(7):659--670, 1976.

\end{thebibliography}
\end{document}